\documentclass[leqno,12pt]{amsart}
\usepackage{amssymb,amsmath}
\newtheorem{theorem}{Theorem}
\newtheorem{proposition}[theorem]{Proposition}
\newtheorem{lemma}[theorem]{Lemma}
\newtheorem{definition}[theorem]{Definition}

\newcommand{\C}{{\mathbb C}}
\newcommand{\R}{{\mathbb R}}
\newcommand{\Z}{{\mathbb Z}}

\newcommand{\Q}{{\mathbb Q}}
\newcommand{\A}{{\mathbb A}}
\newcommand{\CH}{{\mathcal H}}
\newcommand{\CV}{{\mathcal V}}
\newcommand{\CL}{{\mathcal L}}

\newcommand{\CF}{{\mathcal F}}
\newcommand{\CM}{{\mathcal M}}
\newcommand{\CG}{{\mathcal G}}

\newcommand{\Tr}{{\rm Trace}}
\newcommand{\Trace
}{{\rm Trace}}

\author{ by Mich{\`e}le Vergne}
\address{Centre de Math{\'e}matiques Laurent Schwartz, 91128 Palaiseau,
France;  Institut de Math{\'e}matiques de Jussieu, Th{\'e}orie des
Groupes, Case 7012, 2 Place Jussieu, 75251 Paris Cedex 05, France}
\email{vergne@math.polytechnique.fr, vergne@math.jussieu.fr}

\title[All what I wanted to know and was afraid to ask]{All what I wanted to know about Langlands program and
was afraid to ask.}
\date{October 2005}
\begin{document}
\maketitle

\tableofcontents

\newpage

\section{Sources and short description}

My notes are extracted from  informal conversations with G{\'e}rard
Laumon, Eric Vasserot and Jean-Louis Waldspurger and from the
following articles.

Stephen Gelbart: {\it An elementary introduction to the Langlands
program.}  Bulletin of the American Mathematical Society, 10,
1984, pp 177-215.

Several articles (E. Kowalski, Gelbart, Gaitsgory)  from the
book:{ \it An introduction to the Langlands program}. Editors:
Bernstein-Gelbart; Birkhauser (2003).

Thomas Hales: {\it A statement of the fundamental lemma.} Arkiv
RT/031227 (2003).

Gerard Laumon {\it Travaux de Frenkel,Gaitsgory et Vilonen sur la
correspondence de Drinfeld-Langlands.} S{\'e}minaire Bourbaki Juin
2002

 I refer to these articles to get more complete bibliography on
 this subject and more details.

\bigskip

Here is the description of these notes. I will  try  to motivate
the ``Langlands program" by discussing first
 $L$-functions attached
to representations of the Galois group of number fields. Then I
will discuss the Langlands program and its conjectural
implications of functoriality. I will discuss a simple example of
the fundamental lemma, as a combinatorics problem of counting
lattices. Then I will outline the Langlands-Drinfeld program for
the function field of a complex curve. Here this is very very
sketchy.

 I give here only the
definitions I understood. So when I write: ``it is more
complicated", it means I do not understand. I oversimplified many
definitions and conjectures, some willingly, and probably many
unwillingly. Furthermore, I will only discuss here representations
in characteristic $0$, so that I will not touch upon  the recent
developments on modular representations, modular Serre conjecture,
etc, \ldots.

\bigskip

Section \ref{importance} describe some historical motivations to
the Langlands program.

Section \ref{automorphic} gives the main definitions.

The last two sections  ``Functoriality" and ``Geometric Langlands
correspondence" are independent.
  Here I have only be trying to
give some statements and some simple examples. No indications of
techniques  are given. Thus probably, the beauty of the works of
researchers in this field will not  be fully apparent, as the
interplay between representation theory and algebraic geometry is
most remarkable,
 but I felt
too incompetent on these domains (Hitchin moduli spaces, perverse
sheaves, etc..) to write something not totally nonsensical.

\newpage

\section{The importance of automorphic forms}\label{importance}

\subsection{Classical automorphic forms}
 What is a (classical) automorphic form ?  Roughly speaking, this is a function
$f(z)$ on the upper-half plane $\CH:=\{z=x+iy; y>0\}$ of the form
$$f(z):=\sum_{n\geq 0} a_n e^{2i\pi nz/N}$$ ($N$ a positive
integer) and such that it takes almost the same form (I suppose
this is why it has the name automorphic) when changing $z$ to
$-1/z$:
$$f(-1/z)={\rm constant}\cdot z^k f(z),$$ for some $k$.

 Remark that $f$ is periodic of period $N$:
$f(z+N)=f(z)$. The integer  $k$ will be called the weight ($k$ can
be a half integer, but we will mostly restrict ourselves to
integer weights).  The function $|e^{2i\pi nz/N}|=e^{-2\pi ny/N} $
decreases very rapidly as a function of $n$ when  $y={\rm
Im}(z)>0$
 and $n>0$. Provided the coefficients $a_n$ have a
reasonable growth (for example polynomial in $n$), the series is
indeed convergent

For example, the Theta series $\Theta(z)=\sum_{n\in \Z}e^{i\pi n^2
z}$ satisfies $\Theta(-1/z)=  (-iz)^{1/2} \Theta(z)$, as follows
from Poisson formula. Thus it is an automorphic form of weight
$1/2$.

\bigskip

Let us give  more precise definitions. Let $G=SL(2,\R)$ be the
group of holomorphic transformations of the upper-half plane
$\CH$:
 $$G:=\{\left(%
\begin{array}{cc}
  a& b \\
  c & d \\
\end{array}%
\right);  (a,b,c,d) \mbox{ reals}, ad-bc=1\}.$$ The corresponding
transformation is $g\cdot z=(az+b)/(cz+d)$.

Let $N$ be a positive integer, let $\Z/N\Z$ be the ring of
integers mod $N$ and $(\Z/N\Z)^*$ the multiplicative group of
invertible elements of this ring. Let $\chi$ a character of
$(Z/NZ)^*$ and $k$ an integer or half-integer. We extend $\chi$ to
a periodic function of period $N$ of $\Z$ by setting $\chi(u)=0$
if $u$ is not relatively prime to $N$ while
$\chi(u)=\chi(mod(u,N))$ otherwise, so that
$\chi(nm)=\chi(n)\chi(m)$. Such a $\chi$ is called a primitive
Dirichlet character of level $N$  if $N$ is the smallest period.
Most of the time in this introduction, we can think that $N=1$ and
$\chi=1$ (integers prime to $1$ form the empty set, thus $\chi=1$
(Hum!); anyway this is the convention).

 Consider the discrete subgroup of
$SL(2,\R)$:
$$\Gamma_0(N):=\{\left(%
\begin{array}{cc}
  a& b \\
  c & d \\
\end{array}%
\right); (a,b,c,d)\, \mbox{integers}, ad-bc=1,\, c\equiv 0\,
\mbox{mod} N\}.$$
Remark that $\Gamma_0(N)$ contains the matrix
$\left(%
\begin{array}{cc}
  1& 1 \\
  0 & 1 \\
\end{array}%
\right)$ producing the transformation $z\mapsto z+1$ of $\CH$.

If $N=1$, the group $\Gamma_0(1)$ is simply $SL(2,\Z)$.

\begin{definition}
Let $k$ be an integer, $N$ an integer and $\chi$ a character.
 The
space $M(N,\chi,k)$ is (with some analytic restrictions) the space
of functions $f(z)$ of $z\in \CH$ of the form
$$f(z):=\sum_{n\geq 0} a_n e^{2i\pi nz}$$
and such that:
$$f((az+b)/(cz+d))=\chi(a)(cz+d)^kf(z)$$
for all $g=(\begin{array}{cc}
  a& b \\
  c & d \\
\end{array})$ in $\Gamma_0(N)$ .
The integer $N$ is called the level, $k$ the weight.
\end{definition}
Remark that we decided that our  automorphic form $f$ should be
periodic of period $1$: with respect to the preceding informal
definition, where $f$ was periodic of period $N$, we made the
change of variable $z\mapsto z/N$. It is customary to write
$q=e^{2i\pi z}$, so that

$$f=\sum_{n\geq 0} a_n q^n.$$

 If $a_0=0$,
then $f$ is said to be cuspidal at $\infty$: it vanishes on the
cusp $\infty$ of the domain $\Gamma_0(N)\backslash\CH$.

An example is $\Delta\in M(1,1,12)$ with
$$\Delta:=q\prod_{n\geq 1}(1-q^n)^{24}.$$

The space $M(N,\chi,k)$ is a finite  dimensional vector space and
its dimension is known. For example, $M(1,1,12)=\C \Delta$ and
$M(1,1,k)=0$ if $k<12$.

To a cusp form  $f=\sum_{n>0} a_n e^{2i\pi nz}$, one associates
its $L$ series
$$L(s,f)=\sum_{n>0} a_n n^{-s}=\frac{(2\pi)^s}{\Gamma(s)}\int_{0}^{\infty} f(iy) y^{s}\frac{dy}{y}.$$

It is particularly nice to have an Euler product formula for
$L(s,f)$ similar to the formula:
$$\zeta(s)=\prod_{p;  {\rm primes}}\frac{1}{1-\frac{1}{p^s}}$$ for the Riemann
zeta function $\zeta(s)=\sum_{n\geq 1} \frac{1}{n^s}$.

\begin{definition}
Let $\chi$ be a Dirichlet character.
 We say that $f=\sum_{n>0}a_n
e^{2i\pi nz}$ is a newform of weight $k$ and level $N$  if:

1) $$L(s,f)=\prod_{p; {\rm primes}} \frac{1}{1-a_p p^{-s}+\chi(p)
p^{k-1} p^{-2s}}$$

and

2) $z^kf(-1/Nz)$ is proportional to ${\tilde f}(z)=\sum_{n>0}
\overline{a_n} e^{2i\pi nz}$.

We denote by $S(N,k,\chi)$ the space of such newforms.

\end{definition}
 Note that
the coefficients $a_n$ of the expansion
$$f(z)=\sum_{n\geq 1}a_n e^{2i\pi nz}$$ are entirely
determined by the $a_p$ for $p$ prime and $a_1$.

I will ignore the description of the coefficient of
proportionality in condition 2) above, although  it has a great
importance in the theory as
 the ``epsilon" factor.

\bigskip

Using eventually lower levels, Hecke showed (with the help of
operators later
 called Hecke operators, generating the so-called Hecke algebra)
that every cusp form is a linear combination of newforms.

Recall the  definition of Hecke operators.
 Let $p$ be a
prime number and consider the set $D(p)=\{\gamma\in GL(2,\Z);
\det(\gamma)=p\}$.  It is easy that every integral $2$ by $2$
matrix $g$ with integral coefficients of determinant $p$ can be
written  as $g_1 (\begin{array}{cc}
                  p & 0 \\
                  0 & 1 \\
                \end{array}) g_2$, where $g_1,g_2$ are integral matrices of  determinant $1$.
                In other words  $D(p)$ is a double coset for
the action of $SL(2,\Z)$ by left and right translation.
$$D(p)=SL(2,\Z)(\begin{array}{cc}
                  p & 0 \\
                  0 & 1 \\
                \end{array}) SL(2,\Z).$$

It  $f$ is an automorphic form of even weight $k$ (for
$SL(2,\Z)$), then

$$(H_p
f)(z)=p^{-1}\left(\sum_{u=0}^{p-1}f(\frac{z+u}{p})+p^{k}f(pz)\right)$$
is again an automorphic form of weight $k$  and the operators
$H_p$ commute for all primes $p$.

The definition of $H_p$ is  natural when interpreting  an
automorphic form as a function on $SL(2,\R)$. If   $f$ is an
automorphic form of even weight $k$ (for $SL(2,\Z)$), the function
$$(Af)(g)=(a-ci)^{-k}f(g^{-1}\cdot i)=(a-ci)^{-k}f((di-b)/(-ci+b))$$
is a function on $SL(2,\R)$ invariant by right translation by
$SL(2,\Z)$.

Consider the double coset ${\tilde D}(p)$
$${\tilde D}(p)=SL(2,\Z)\left(\begin{array}{cc}
                  p^{1/2} & 0 \\
                  0 & p^{-1/2} \\
                \end{array}\right) SL(2,\Z)$$
in $SL(2,\R)$. As ${\tilde D}(p)$ is left and right invariant by
$SL(2,\Z)$, if $F$ is a function of $SL(2,\R)$  invariant by right
translation by $SL(2,\Z)$, then
$$({\CH_p}F)(g)=\frac{1}{p}\sum_{\gamma\in {\tilde
D}(p)/SL(2,\Z)}F(g\gamma)$$ is still right invariant by
$SL(2,\Z)$.

The equality  $$p^{k/2}\CH_p Af=A(H_pf)$$ provides a natural
interpretation of $H_p$ as  summing a function over the finite
$SL(2,\Z)$-orbit ${\tilde D}(p)/SL(2,\Z)$ in $SL(2,\R)/SL(2,\Z)$.

Let us check this equality. Let $F=Af$ and compute $
F(g(x,y)\gamma(u))$ with $$g(x,y)=\left(\begin{array}{cc}
                  y^{-1/2} & -xy^{-1/2} \\
                  0 & y^{1/2} \\
                \end{array}\right),\hspace{1cm}
                \gamma(u)=\left(\begin{array}{cc}
                  1 & -u \\
                  0 & 1 \\
                \end{array}\right)
\left(\begin{array}{cc}
                  p^{1/2} & 0 \\
                  0 & p^{-1/2} \\
                \end{array}\right).$$

We obtain
$$p^{k/2}F(g(x,y)\gamma(u))=y^{k/2}f(\frac{z+u}{p}).$$

If $$\gamma_0= \left(\begin{array}{cc}
                  p^{-1/2} & 0 \\
                  0 & p^{1/2} \\
                \end{array}\right),$$

we obtain  $$p^{k/2}F(g(x,y)\gamma_0)=p^{k}y^{k/2}f(pz).$$

   The equality $p^{k/2}\CH_p Af=A(H_pf)$ follows from the fact
   that the elements $\gamma(u), u=0,\ldots, (p-1)$ and $\gamma_0$ are
   the representatives of the finite set ${\tilde
D}(p)/SL(2,\Z)$.

\bigskip

A newform is an eigenvector for all operators $H_p$. Furthermore,
the Fourier coefficient $a_p$ of a newform (with $a_1=1$) is
exactly the eigenvalue of $H_p$. We record this important fact in
a proposition.

\begin{proposition}\label{Hecke}
If $f$ is a newform with $a_1=1$, then
 $$H_pf=a_p f.$$
\end{proposition}

The automorphic property of $f\in S(N,k,\chi)$ translates
immediately to a relation between $L(s,f)$ and $L(k-s,\tilde{f})$:
``proof": change $y$ in $1/Ny$ in the integration formula.  This
gives for example the proof of the functional equation for the
Riemann Zeta function from  the $\Theta$ series automorphic
property. More precisely, define for a newform $f$,
$$\Lambda(s,f)=(2\pi)^{-s}\Gamma(s)L(s,f).$$
Then one obtain
$$\Lambda(s,f)=\epsilon(f)N^{k/2-s}\Lambda(k-s,{\tilde f})$$
with a constant factor $\epsilon(f)$, the ``$\epsilon$" factor.

The hope  of  the ``automorphic  gang" is that any function $L$
(defined as a nice looking product over primes) arising in the
world of mathematics is the $L$ function of an automorphic form,
in short $L$ is modular. This optimistic assumption had already
many successes. Recall for example that the word ``END" appears on
Wiles proof of Fermat's theorem after the phrase: If there is a
solution of Fermat equation, then there is a  newform of weight
$2$ and level $2$   but there is no such forms. END.

\subsection{ Abelian reciprocity law}\label{abelian}

\subsubsection{A few notations}\label{notation}

Let $G$ be a group.

Two elements $g$ and $g'$ are conjugated in $G$ if there exists
$u\in G$ such that $g'=ugu^{-1}$. If $G$ is abelian (commutative),
then $g'=g$. Otherwise, the conjugacy class $O_g$ of an element
$g$ in a group $G$ is the subset $$O_g:=\{ugu^{-1}, u\in G\}$$ of
all the conjugated elements to $g$ in $G$.

A representation $\sigma$ of a group $G$ (finite or infinite)
associates to $g\in G$ a linear transformation $\sigma(g)$ of a
vector space $V$. We must have $\sigma(1)=Id_V$ and $\sigma(g_1
g_2)=\sigma(g_1)\sigma(g_2)$. In this report, our spaces $V$ will
be complex vector spaces.

 If $V$ is one dimensional, then $\sigma$ is a character of $G$.
In general, $\sigma(g)$ is represented by a  matrix in a basis of
$V$. We can take (at least when $V$ is finite dimensional) its
trace ${\rm Trace}_V(\sigma(g))$ and its characteristic polynomial
$\det_V(1-\sigma(g))$. Remark that the trace and characteristic
polynomial take same value on two conjugated elements $g,g'$.

Most of the time, if $V$ is a vector space provided with an action
of $G$, $v$ a vector of $V$ and $g$ an element of $G$, we note
$\sigma(g)\cdot v$ simply by $g\cdot v$.

If $V$ is a representation of $G$ and $H$ a subgroup of $G$, we
define the space  $V^H$ of $H$ invariants by

$$V^H:=\{v\in V; h\cdot v=v \,{\rm for\, all}\, h\in H\}.$$

\bigskip

Let $p$  be a prime  and $F$ a finite extension of the field
$\Z/\Z p$. Then $x\mapsto x^p$ is an automorphism of $F$, called
the Frobenius ${\rm Fr }$: the fixed elements of ${\rm Fr}$ is the
ground field $\Z/p\Z$ ($a^p=a \,{\rm mod}\, p$ for
 all $a\in \Z$).

\subsubsection{ Abelian reciprocity law}

The motivations of Langlands program  come from number theory.

Let $E$ be a number field: a finite extension of $\Q$. We consider
the ring $O_E$ of integers of $E$: that is the set of solutions
$x$ of an equation $x^n+\sum_{i<n} a_i x^i=0$ with $a_i\in \Z$. A
prime $p\in \Z$ gives rise to an ideal $(p)=O_Ep$ of $O_E$. Assume
$E$ is a Galois extension, with Galois group $G$. The ideal $(p)$
factors as $(\prod_i P_i)^e$ where $P_i$ are different prime
ideals in $O_E$ and $e$ is called the ramification index of $p$.
If $e=1$, then $p$ is said to be unramified. The prime ideals
$P_i$ are said to be above $p$. Assume $p$ is unramified (all $p$,
but a finite number), and choose $P$ above $p$. Then there exists
a unique  element ${\rm Fr}_P$ of the Galois group $G$ such that

\begin{itemize}
\item

The element ${\rm Fr}_P$ leaves stable the prime $P$ above $p$.\\
\item

${\rm Fr}_P$ induces the Frobenius transformation
in the field  $O_E/P$ extension of $\Z/p\Z$.\\
\end{itemize}

For example, if $p$ splits completely (number of $P$ above
$p$=degree of the extension), then  ${\rm Fr}_P=I$.

For $p$ unramified, when the prime $P$ varies above $p$,  the
Frobenius element ${\rm Fr}_P\in G$ are in the same  conjugacy
class, and sometimes we will write ${\rm Fr}_p$ for an element in
the conjugacy class of ${\rm Fr}_P$. If $G$ is abelian, ${\rm
Fr}_p$ is well defined.

If $p$ is ramified, choose a prime $P$ above $p$, let $D_P$ the
stabilizer of the prime ideal $P$ in the Galois group $G$. The
group $D_P$ is called the decomposition group. There is a morphism
of $D_P$ onto the Galois group of $O_E/P$ with kernel $I_P$, the
inertia group. The cardinality of $I_P$ is the ramification index
$e$. Now ${\rm Fr}_P$ is still defined as the element producing
the Frobenius transformation in the field  $O_E/P$ extension of
$\Z/p\Z$, but now it leaves in $D_P/I_P$.

-{\bf Example}
 {\it Take the example $E=Q[\sqrt{-1}]=Q[i]$ with ring of integers
 $O_E=\Z+\Z i$.
The Galois group is $G=\Z/2\Z$, generated by the complex
conjugation $i\mapsto -i$. If $p\neq 2$ is the sum of $2$ squares
$p=n^2+m^2$, then $(p)=(n+mi)(n-mi)$, thus $p$ splits completely
and ${\rm Fr}_p=I.$

If $p\neq 2$ is not a sum of two squares, then the extension
$O_E/O_Ep$ is an extension of degree $2$ and  ${\rm
Fr}_p=Conjugation$.

We have $(2)=(1+i)(-i(1+i))$, so that $(2)=(\Z[i](1+i))^2$, and
$2$ is the only ramified prime.

In this example, ${\rm Fr}_p=I$ is equivalent to
 $p=1$ mod $4$.}

\bigskip

It is very important to understand ${\rm Fr}_P$ as we want to
understand the ``factorization" of primes.

The abelian reciprocity law is the following.

\begin{theorem} {\bf Abelian reciprocity law} (KRONECKER-WEBER).

Let $E$ be a Galois extension of $\Q$ with Galois group $G$. Let
$\sigma: G\to \C^*$ a character of the Galois group. Then there
exists an integer $N_\sigma$ (the conductor) and a primitive
Dirichlet character $\chi_\sigma$ of $\Z$ of level $N_\sigma$ such
that
$$\sigma({\rm Fr}_p)=\chi_\sigma(p)$$
for all unramified $p$.

\end{theorem}

{\bf Example}

{\it Take $\Q(i)$ with $G= \{1,conjugation\}$ as Galois group.
Thus for the character $\sigma:G\mapsto \{1,-1\}\subset \C^*$,
sending complex conjugation to $(-1)\in \C^*$, we have
$\sigma({\rm Fr}_p)=\chi_\sigma(p)$ where $\chi_\sigma$ is the
primitive Dirichlet character mod $4$ defined by
$\chi_\sigma(n)=(-1)^{(n-1)/2}$, when $n$ is odd.}

\bigskip

Let $E$ be an abelian extension of $\Q$ and $\sigma$ a character
of $G$. To each prime  integer $p$,  introduce the $L$-factor
$L_p(s,\sigma)$. When $p$ is unramified:
$$L_p(s,\sigma)= \frac{1}{(1-\sigma({\rm Fr}_p)p^{-s})} .$$
For ramified $p$, define $$L_p(s,\sigma)=1.$$ Then define the
$L$-function
$$L(s,\sigma)=\prod_{p; \,{\rm primes}\,}L_p(s,\sigma).$$

From the Abelian reciprocity law, we see that

$$L(s,\sigma)=\sum_{n\geq 1} \chi_\sigma(n) n^{-s}$$
and again Poisson formula shows that $L(s,\sigma)$ is entire, if
$\chi\neq 1$, and has some functional equation.

To summarize:

The abelian reciprocity law relates the Galois group of an abelian
extension to the groups $(\Z/N\Z)^*$  which are the Galois group
of the cyclotomic field  $\Q(\mu_N)$ obtained by adding to $\Q$ a
primitive $N$-th root of unity $\mu_N$. It has as consequence that
any abelian extension can be imbedded in a cyclotomic field:

{\it Example: $\Q(\sqrt{2})\subset Q(\mu_8)$ as
$\sqrt{2}=e^{i\pi/4}+e^{-i\pi/4}$.}

\subsection{Artin $L$-functions and  Langlands conjectural
reciprocity law}\label{langlands}

\subsubsection{Artin L functions} Langlands aim was to formulate
and eventually prove a ``non-abelian" reciprocity law.

Now consider representations $\sigma$  of $G$, where $G$ is the
Galois group of a Galois extension of $\Q$ in a finite dimensional
complex vector space $V$ of dimension $n$. Thus each element $g\in
G$ is represented by a $n\times n$ complex matrix.

To each prime  integer $p$,  introduce the local $L$-factor
$L_p(s,\sigma)$, as follows. Take a prime $P$ above $p$.

If $p$ is unramified, define

$$ L_p(s,\sigma)=\frac{1}{\det(1-\sigma({\rm Fr}_P) p^{-s})}.$$

If $p$ is ramified, we consider the space $V^{I_P}$ of invariants
of $V$  under the inertia group $I_P$. Then we can define
$$L_p(s,\sigma)=\frac{1}{\det(1-\sigma({\rm Fr}_P) p^{-s})}$$
where the transformation ${\rm Fr}_P$ acts on the space of $I_P$
invariants.

The definition of the local $L$-factor does not depend of the
choice of $P$ above $p$.

Artin defined

\begin{equation}\label{L}
L(s,\sigma)=\prod_{p; \,{\rm primes}\,} L_p(s,\sigma).
\end{equation}

\bigskip

If $\sigma$ is of dimension $2$, and $p$ is unramified, the
element ${\rm Fr}_p$ (defined only up to conjugacy) of $G$ is send
to a matrix $\sigma({\rm Fr}_p)$ with two eigenvalues $\alpha_p,
\beta_p$. The local $L$ factor associated by Artin to $\sigma$  is
thus $$ L_p(s,\sigma)=\frac{1}{\det(1-\sigma({\rm Fr}_p) p^{-s})}=
\frac{1}{1-\Tr(\sigma({\rm Fr}_p)) p^{-s}+\det(\sigma({\rm
Fr}_p))p^{-2s}}.$$

 Remark that $g\mapsto \det(\sigma(g))$ is a one dimensional
character of $G$, so that by the Abelian reciprocity theorem, we
already have the existence of a Dirichlet character $\chi$  of
level $N$ such that
$$\det(\sigma({\rm Fr}_p))=\chi(p).$$

Thus $$L_p(s,\sigma)=\frac{1}{1-(\alpha_p+\beta_p)
p^{-s}+\chi(p)p^{-2s}}.$$

Artin  conjectured that $L(s,\sigma)$ is entire, when $\sigma$ is
irreducible.

Let $\sigma$ be an odd representation of $G$ (help ?).
 The conjecture of  Artin (still unproved) can be restated by
saying that  $L(s,\sigma)$
 is automorphic: more precisely there exists
a newform $f\in S(N,\chi,1)$ of level $N$, weight one and
Dirichlet Character $\chi$ such that
$$L(s,\sigma)= L(s,f).$$

In short, the traces  ${\rm Trace} (\sigma({\rm Fr}_p))$ of the
Frobenius elements should be the Fourier coefficients $a_p$ of a
newform $f=\sum_{n\geq 1} a_n e^{2i\pi nz}$ for all unramified
prime $p$.

%
%
%
%

\subsubsection{Langlands reciprocity conjecture}

Let now $\sigma$ a representation of the Galois group in $\C^n$.
 Langlands formulated the conjecture that  $L(s,\sigma)$ is the
$L$-function associated to an automorphic representation of
$GL(n,\A)$ where $\A$ is the ring of adeles of $\Q$. {\bf  Roughly
speaking representations $\sigma$ of degree $n$ of the absolute
Galois group $G({\overline \Q}/\Q)$ parametrize some automorphic
representations of $GL(n,\A)$}. We give more details later.

 Langlands-Tunnell (\cite{langlands},\cite{tunell})  have proved that the Artin
$L$-function for representations of the Galois group on $\C^2$
with image a solvable group are $L$-functions of automorphic
forms. Already to settle  this case, base change, the trace
formula, the fundamental lemma, the lifting of automorphic
representations have to be established for non trivial cases. The
Langlands program has taken a life of its own since then, and many
results have been proved, as parts of the Langlands original
``program" or inspired by it. Indeed, many natural problems:
functoriality, base change, local correspondance arise from this
dictionary (representations of the Galois group = automorphic
forms). It is clear that now the study of automorphic forms is a
central topic in mathematics, with interconnections with algebraic
geometry, arithmetic geometry, representations of quantum groups,
etc...

Up to now, we discussed the field $\Q$. It is also important to
have the same theory for any number field. For example, Langlands
proof of the Artin's conjecture for the solvable subgroup $A_4$
($A_4\subset PGL(2,\C)=SO(3,\C)$ embedded as the symmetry group of
the tetrahedron) uses a composition series leading to study cubic
extensions $E$ of $\Q$ and the corresponding base change.

So let $E$ be a number field. Similarly  a representation of the
Galois group $G(\overline{E}/E)$ of dimension $n$ should lead to
an automorphic representation of $GL(n,\A_E)$, where now $\A_E$ is
defined using all completions of $E$. If $E$ is a number field,
then clearly  $G(\overline{\Q}/\Q)$ has an homomorphism into the
finite group  $G(E/\Q)$ with kernel $G(\overline{\Q}/E)$. Recall
that if $K$ is a given finite group, it is unknown if $K$ is a
Galois group of a number field. So  the idea to study
representations of $G(\overline{\Q}/\Q)$ with finite image may be
a tool to understand the possible   Galois groups.

\subsubsection{Arithmetic varieties and Langlands
conjecture}\label{aritm} (very sketchy)

There are many representations of the Galois group $G=G(\overline
\Q/\Q)$  occuring in ``nature".

Let $X$ be an arithmetic variety: the set of solutions of
equations defined over $\Z$. There are very natural
representations of the absolute Galois group $G$ associated to
$\ell$-adic cohomology groups of  $X$ (these representations do
not factor through finite groups and do not have finite images in
$GL(n,\C)$), but with the same formulae (\ref{L}) as in the
preceding section, they give rise to $L$-functions (called motivic
$L$-functions). Furthermore if $X$ is provided with an action of a
group $S$, then the $\ell$-adic cohomology groups $H_{\ell}^i$ are
provided with a representation of $S\times G$.

  For
example, let $X:=\{y^2=4x^3-Ax-B, A,B\in \Z,A^3-27 B^2\neq 0 \}$ a
smooth elliptic curve, then its first $\ell$-adic cohomology group
is a vector space of dimension $2$ over $\overline{\Q_{\ell}}$ and
this representation of dimension $2$ of $G$ gives rise to the
$L$-function attached by Hasse-Weil to $E$, described
``concretely" as follows: we consider the number of points $N(p)$
of $X$ in the finite field $\Z/p\Z$. Then for good primes $p$, the
local $L$ factor is $\frac{1}{1-a_pp^{-s}+p p^{-2s}}$, with
$N(p)=1-a_p+p$. The Shimura-Taniyama-Weil conjecture, proved by
Wiles and Taylor (+ Diamond, Conrad, Breuil), says that these
$L$-functions are modular.

Another example. Let $X=X_0(N)$ be the Borel-Bailey
compactification of $\Gamma_0(N)\backslash\CH$, then $X$ is an
arithmetic curve of genus $g$, called the modular curve. Thus
there are exactly $g$ newforms of weight $2$ and character
$\chi=1$.
 The $L$-function attached to
the first $\ell$-adic cohomology group $H^1_{\ell}$ (a vector
space of dimension $2g$) is the product of the $L$ functions
$L(s,f_i)$ where $f_1,f_2,\ldots,f_g$ are the $g$ newforms. Here
for each prime $p$, the Hecke operator $H_p$  as well as the
Frobenius element $Fr_p$  acts on $H_{\ell}^1$ and have the same
eigenvalues.  This gives ``the" explanation of the correspondence
between automorphic forms and some representations of $G$.

Higher dimensional analogues of the modular curves are Shimura
varieties (for the symplectic group ${\rm Sp}(n,\Z)\backslash {\rm
Sp}(n,\R)/U(n)$). Unfortunately for $GL(n)$ and $n>2$, there is no
analogue of the Shimura varieties.

\section{Automorphic representations}\label{automorphic}

\subsection{The use of adeles}

Interpreting the $L$ function of Artin associated to
representations $\sigma$ of $G$ in $GL(n,\C)$ needs the notion of
automorphic representation. If $n=2$, there is a nice theory
(classical) of automorphic forms as  part of the function theory
of the upper-half plane. But to generalize it to $GL(n,\R)$, it is
easier to use adeles. The factorisation of the $L$ function over
primes will then have natural interpretations, etc...

If $p$ is a prime number, a $p$-adic integer is a
 series $\alpha=\sum_{n=0}^{\infty} a_p p^n$ with $a_p$ an integer
between $(0,p-1)$. Thus $$\Z_p:=\{\alpha=\sum_{i=0}^{\infty} a_i
p^i, 0\leq a_i<p\}$$ and the ring of fractions of $\Z_p$ is
naturally identified to
$$\Q_p=\{\alpha=\sum_{i>i_0}^{\infty} a_i p^i, 0\leq a_i<p\}$$
where now $i_0$ can be negative.

As $\Z$ is naturally embedded in $\Z_p$ by writing an integer in
base $p$, the ring $\Q$ is embedded in $\Q_p$.

 The order of an element  $\alpha$  in $\Q_p$, $\alpha$ non zero,
  is the smallest $i$ with $a_i\neq 0$ and is
 denoted by $val(\alpha)$.
 The integer $p$ has order $1$, and called the uniformizer.
The group  $\Z_p^*$ is exactly the set of elements of order $0$ in
$\Q_p$ and $\Q_p^*/\Z_p^*=\Z$, via the map $val$.

The topology on $\Q_p$ is as follows. A sequence $x_n$ converges
to $0$, if the orders of the elements $x_n$ tends to $\infty$
(more and more divisible by $p$). Thus it is clear that $\Z_p$ is
compact and open. The normalized  additive Haar measure of $\Q_p$
gives mass $1$ to $\Z_p$. Thus on $\Q_p/\Z_p$, integrating means
counting. For example, when $r>0$,  the measure of the set $\{u,
val(u)\geq -r\}$ modulo $\Z_p$ is $p^r$. The Haar measure on the
locally compact multiplicative group $\Q_p^*$ is normalized by
giving mass $1$ to  $\Z_p^*$. Thus integration on
$Q_p^*/\Z_p^*=\Z$ means counting.

The adele ring $\A$ is the ring $\R\times\prod_{\{p;  \,{\rm
primes}\,\}} \Q_p$ where we assume that  almost all components
$\alpha_p$ of an element $\alpha=(\alpha_p)$ are in $\Z_p$.

Then  $\Q$ is embedded diagonally in $A$:
$$\alpha\mapsto (\alpha,\alpha,\alpha,\ldots)$$

Let $\A^*$ be the group of invertible elements of $\A$. Clearly,
we can multiply an element $a\in \A^*$ by some element $\alpha\in
\Q^*$ so that $\alpha a$ is  in $\R^+\times \prod_p \Z_p^*$. Thus
$$\A^*=\Q^*\times \R^+\times \prod_p \Z_p^*.$$

\subsection{A little more of representation theory}

If $G$ is a (locally compact) group, a very natural representation
of $G$ is the regular representation of $G$. It acts on the space
$L^2(G)$ of  functions (square integrable) on $G$ by
$$R(g)f(u)=f(g^{-1}u).$$

In particular if $f$ is a function in $L^2(G)$, the closed suspace
${\rm Translate}(f)$ generated by linear combinations of
translates of the function $f$ is a $G$-invariant subspace of
$L^2(G)$.

More generally, if $\Gamma$ is a subgroup of $G$, the space
$L^2(G/\Gamma)$ of functions on $G$ such that $f(g\gamma)=f(g)$
for all $g\in G$, $\gamma\in H$ is the ``most natural" way to
construct representations of $G$. If $f$ is in $L^2(G/\Gamma)$,
the space ${\rm Translate}(f)$ is a subspace of $L^2(G/\Gamma)$.
The group $G$ acts on $L^2(G/\Gamma)$ by left translations
$L(g_0)f(g)=f(g_0^{-1}g)$. The corresponding representation of $G$
in $L^2(G/\Gamma)$ is called a quasi regular representation ( or a
permutation representation) and is denoted by $Ind_{\Gamma}^G 1$.

A representation $\pi$ of $G$ in a Hilbert space  $V$ is
irreducible, if $V$ does not admit (closed) non trivial invariant
subspaces.

\subsection{From classical automorphic forms to automorphic
representations}\label{aut}

Here we consider first $GL(2)$. We write $\CV:=\{p; {\rm prime\,
integers}\}\cup \{\infty\}$, the set of valuations. The local
groups are  $GL(2,\Q_p)$  for $p$ a prime and $GL(2,\R)$ for
$v=\infty$.

 We write $K_p:= GL(2,\Z_p)$ (in particular  $\det(g)\in \Z_p^*$, for $g\in GL(2,\Z_p)$).
An element $g$ of the group $GL(2,\A)$ is a family $(g_v)_{v\in
\CV}$ where for all prime $p$, except a finite number, $g_p$ is in
$K_p=GL(2,\Z_p)$. Similarly we can send $GL(2,\Q)$ in $GL(2,\A)$
where it becomes a discrete subgroup.
 We write $K_0=\prod_{p<\infty}K_p$. The center of $GL(2,\A)$ is $$Z_\A=\{\left(%
\begin{array}{cc}
  a& 0 \\
  0 & a \\
\end{array}%
\right);  a\in \A^*\}.$$
 Then
$$K_0\backslash\ GL(2,\A)/GL(2,\Q)Z_\A= SL(2,\R)/SL(2,\Z).$$

If $f$ is a classical  automorphic form on $\CH$ (for $SL(2,\Z)$)
of weight $k$,  it is easy to see that
$\phi(g)=f(g^{-1}.i)(ci+d)^{k}$ is a function on $
SL(2,\R)/SL(2,\Z)$. Here $g=(\begin{array}{cc}
  a& b \\
  c & d \\
\end{array})$.
Then given a classical newform $f$ (for $SL(2,\Z)$) on $\CH$ of
weight $k$, there exists a unique function $\CF$ in $L^2(
GL(2,\A)/GL(2,\Q)Z_\A)$ such that  $\CF$ coincide with $\phi(g)$
on $SL(2,\R)$ and such that $\CF$ is invariant by left
translations by all subgroups $K_p$ for primes $p$. Then the
function $\CF$ is an element of $L^2(GL(2,\A)/GL(2,\Q))$ and the
space ${\rm Translate}(\CF)$  span an irreducible subspace in
$L^2(GL(2,A)/GL(2,\Q)Z_\A)$: this is by definition an automorphic
representation of $GL(2,\A)$.

We summarize: A newform $f$ (for $SL(2,\Z)$) is the same thing as
a function $\CF$ on $GL(2,\A)/GL(2,\Q)Z_\A$ such that $\CF$ is a
decomposable vector $\CF=\otimes_{v\in \CV}\CF_p$, and $\CF_p$ is
left invariant invariant by all compact groups $GL(2,\Z_p)$ for
all primes $p$.

This definition can easily be generalized to any integer $n$. We
define $GL(n,\A), GL(n,\Q)$ as before, with embedding $
GL(n,\Q)\mapsto GL(n,\A)$.

\begin{definition}
 An automorphic cuspidal representation $\pi$ of $GL(n,\A)$ is an
irreducible subrepresentation of $L^2(GL(n,\A)/GL(n,\Q)Z_\A)$.
Furthermore the representation must satisfy a certain ``cuspidal"
condition.
\end{definition}

In fact, the representation  $\pi$ is necessarily given as a
product $\prod_{v\in \CV} \pi_v$ where $\pi_v$ are irreducible
representations of $GL(n,\Q_v)$. {\bf There is a $L$-function
associated to $\pi$}. Indeed for almost all primes $p$, the
representation $\pi_p$ has a fixed vector under the maximal
compact group $K_p$ and  gives rise to a representation of the
Hecke algebra $\CH_p$ in $\C$, and we will see (next subsection)
that this gives rise to a $L$-factor.

There are similar definitions for number fields of $GL(n,\A_E)$,
etc..

 {\bf General reciprocity law: Langlands conjecture.}

{\bf Let $E$ be a finite extension of $\Q$ with Galois group $G$
and
 $\sigma$ be an irreducible representation of $G$ in $\C^n$.
 Then there exists an automorphic cuspidal  representation $\pi_\sigma$ of
 $GL(n,\A_E)$ such that
 $$L(s,\sigma)=L(s,\pi_\sigma).$$}

Recall that this conjecture is open even for $n=2$ and $E=\Q$.

\subsection{Hecke algebras}

A common tool for studying automorphic representations  is the
Hecke algebra. We already discussed this in the case of classical
automorphic form.

Let $G$ be a locally compact group and  $K$ a compact subgroup.
The algebra $\CH(G,K)$ is the algebra of (compactly supported)
functions on $G$ invariant by left and right translations. This is
an algebra under convolution:
$$(\phi_1*\phi_2)(g)=\int_{g_1,g_2; g_1g_2=g}\phi(g_1)\phi(g_2)$$

Clearly if $(\pi,V)$ is a representation of $G$, the action of the
operators $\pi(\phi)$ (Definition \ref{Tphi}) for $\phi\in
\CH(G,K)$ leaves stable the space $V^K$ of $K$-invariant vectors
in $V$.

Let $G=GL(n,\Q_p)$ and $K=K_p=GL(n,\Z_p)$. Then the algebra
$\CH(G,K)$ is called the Hecke algebra. It has the following
description $\CH(G,K)=\C[T_1,T_2,\ldots,T_n,T_n^{-1}]$ where $T_i$
is the characteristic function of the double coset of the matrix
$$h_p^i=\left(%
\begin{array}{ccccc}
  p & 0 & 0 & 0 & 0 \\
  0 & p & 0 & 0 & 0 \\
  0 & 0 & p & 0 & 0 \\
  0 & 0 & 0 & * & 0 \\
  0 & 0 & 0 & 0 & 1 \\
\end{array}%
\right)$$  with $i$ elements $p$ on the diagonal and the remaining
$(n-p)$ entries equal to $1$.

By definition, an unramified irreducible representation of
$GL(n,\Q_p)$ is a representation having a fixed vector by $K_p$.
The following proposition follows.
\begin{proposition}(Satake)
If $\pi$ is a unramified irreducible representation of
$GL(n,\Q_p)$, then the space of $K$-fixed vectors is
$1$-dimensional and gives rises to a one dimensional $\chi$
character of $\CH(G,K)$.
\end{proposition}

 The local $L$-factor attached to the representation $\pi$ is by
definition

$$L_p(s,\pi)=\prod_{i=1}^n \frac{1}{1-a_ip^{-s}}$$
where  the $i$-th symmetric function $s_i(a)$  of the $a_i$ is
equal to $\chi(h_p^i)$ (up to some power of $p$).

The conclusion is:

Let $\CF=\otimes_v \CF_v$ be a decomposable vector in
$L^2(GL(n,\A)/GL(n,\Q))$ which generates an automorphic
representation. Assume that at $p$, the function $\CF_p$ is left
invariant by $K_p$ ($p$ unramified for almost all $p$), then the
function $\CF_p$  on $GL(n,\Q_p)/GL(n,\Z_p)$ is an eigenvector for
the action of the Hecke algebra $\CH(GL(n,\Q_p)), K_p)$.

\subsection{Local Langlands conjecture.}\label{local}

 Let  $G$ be the absolute Galois group $G({\overline \Q}/\Q)$.
The preceding conjecture says that  there is a map from
representations of $G$ of degree $n$ to automorphic
representations of $GL(n,\A)$. It suggests that at the local
level, irreducible representations $\pi_v$ of $GL(n,\Q_v)$ should
have something to do with the set $Rep_n$ of representations of
degree $n$ of the Galois group $G({\overline \Q_v}/\Q_v)$. In
fact, the right group to take is the Weil   group $W_v$ (strongly
related to the Galois group). This local correspondence has been
proved.

\begin{theorem}{\bf Local Langlands correspondence}
(Harris-Taylor \cite{harris}; simplified by
Henniart\cite{henniart}) Let $W_p$ the local Weil group. There
exists a bijective map
$$Langlands: Rep_n(W_p)\to  Irr(GL(n,\Q_p))$$ such that:
$$L_p(s,\sigma)=L_p(s,Langlands(\sigma))$$
$$(\epsilon_p(s,\sigma)=\epsilon_p(s,Langlands(\sigma)))$$
\end{theorem}

Note that the Langlands conjecture makes sense for any reductive
group $G$. For this, we need  to introduce $L$ groups. I will not
do it here. The local Langlands conjecture for general reductive
groups remains open.

The archimedean local Langlands conjecture was proven by
Langlands.

\subsection{Global Langlands conjecture on function
fields.}\label{lafforgue}

(very sketchy)

 Let $p$ be a prime, and let $F_p=\Z/p\Z$ be the finite field with $p$ elements.
As emphasized by Weil, there is a complete  analogy between number
fields and finite extensions of the field $F_p(t)$. Such an
extensions is the field of rational functions on an algebraic
curve $C$ defined over $F_p$.

\begin{itemize}

\item

\framebox{A Number field $F$} $\mapsto$ \framebox{ A finite \
extension of  $F_p(t)$} $\mapsto$ \framebox{A curve $C$}
\\

\item

\framebox{A prime number} $\mapsto$ \framebox{An irreducible
polynomial in $F_p[T]$}$\mapsto$ \framebox{Points in $C$}
\\

\end{itemize}

The adele ring can be defined, and the theory of automorphic forms
have been developed  for any global field (number fields or
function fields).

The Langlands program can be formulated. In fact, in this case,
there are more geometric tools. In particular, Drinfeld
constructed an arithmetic variety $X$ (chtoucas) over $C\times C$,
where Hecke operators and Frobenius operators acts and could prove
that they have same eigenvalues in the case of  $GL(2)$. Lafforgue
has obtained the proof of the {\bf Global Langlands
correspondence} over function fields in characteristic $p$ for
$GL(n)$.

In the last section, we will discuss very briefly  ``geometric
Langlands correspondence" for a function field of characteristic
$0$.

\section{Functoriality}\label{functoriality}

\subsection{The problem of liftings}

The  Langlands dictionary: ``(Galois representations)=(automorphic
representations)" suggests that some ``trivial operations" in one
side (restrictions of representations) have a counterpart on the
other side (lifting of automorphic forms).

For example, existence of base change and liftings are conjectured
to exists from this correspondence. Let us explain.

{\bf Base Change.} If $E$ is a Galois extension of $\Q$, there is
a map from $G(\overline{\Q}|E)$ to  $G(\overline{\Q}|\Q)$, so that
a representation of $G(\overline{\Q}|\Q)$  gives us a
representation of $G(\overline{\Q}|E)$. Thus there  should be a
base change (denoted $BC$) from automorphic representations $\pi$
of $GL(n,A_\Q)$
  to automorphic representations $BC(\pi)$ of
$GL(n,A_E)$ (and such that the action of $G(E|\Q)$ on $GL(n,A_E)$
leaves fixed the isomorphism class of $BC(\pi)$) .

{\bf Liftings.} Consider a homomorphism $h:GL(n,\C)\to GL(m,\C)$.
For example, if $g\in GL(n,\C)$, then $g$ acts on the space
$S^h(\C^n)$ of homogeneous polynomials  of degree $h$ in $n$
variables by $(g\cdot P)(u)=P(g^{-1}u)$. If $n=2$, and $h=2$, from
a representation of $G$ in $\C^2=\{(x_1e_1+x_2e_2)\}$ we obtain a
representation of $G$ in $\C^3$ (basis $x_1^2,x_1x_2,x_2^2$).
These representations are referred as the symmetric tensors.

If $\sigma$ is a representation of $G$ in $GL(n,\C)$, composing
with $h$, we obtain  a representation of $G$ in $GL(m,\C)$. Thus
according to the  Langlands dictionary, there should be a map
${\rm Lift}_h$ associating to $\pi$ an automorphic representation
of $GL(n,\A)$ an automorphic representation ${\rm Lift}_h(\pi)$ of
$GL(m,\A)$. In other words, there should be a lifting of
automorphic representations, with of course correspondence of the
local factors.

The lifting corresponding to the symmetric tensor $ GL(2,\C)\to
GL(3,\C)$ has been established by Gelbart-Jacquet and
Piateskii-Shapiro. These liftings are highly non trivial.
 If the lifts for any symmetric tensor of a representation of
$GL(2)$  was constructed, as well as general base change for
$GL(2)$, then Artin conjecture would follow. We know a priori what
the lift produces at the level of $L$ functions, but it is quite
difficult to understand what is the lift for representations.  Up
to now, tools are ``special": for example, use of tensor products
of the Weil representations, of the trace formula (with use of the
fundamental lemma), etc...

\subsection{Trace formula}\label{trace}

\subsubsection{Baby version of the trace formula} Let $G$ be a
finite group, and $\phi$ be a function on $G$. Let $g$ be an
element of $G$. The orbital integral  $<O_g,\phi>$ of $\phi$ is
the sum of the values of $\phi$ on the conjugacy class $O_g$ of
$g$:

$$<O_g,\phi>=\sum_{g'\in O_g}\phi(g').$$

 Let $\pi$ be a  finite dimensional representation of $G$.
Then the most important invariant of $\pi$ (this determines $\pi$)
is its ``character": This is the function $g\mapsto {\rm Trace}
(\pi(g))$. Remark that this function is constant over conjugacy
classes, so that the representation $\pi$ is completely determined
by the value of its characters over the set of conjugacy classes.

{\it  Example: the regular representation $R$ of a finite group
$G$: all elements of $g$ except the identity  shifts the elements
of the group, so that $\Trace\, (R(g))=0$ if $g\neq 1$ while
$\Trace (R(1))=|G|$.}

It is useful to introduce the trace as a ``distribution", that is
a linear form on functions on $G$:

$$<{\rm Trace\,}\pi,\phi>=\frac{1}{|G|}\sum_{g\in G} \phi(g){\rm Trace\,}(\pi(g)).$$

{\bf Example.} We have
$${\rm Trace}_{L^2(G)}(g)={\rm Dirac}_1(g).$$

Let $\Gamma$ be a subgroup of $G$ and let us consider the quasi
regular representation $\pi=Ind_\Gamma^G1$. To compute the trace
of this representation on an element $g\in G$, we need to find the
fixed classes $w$: $gw \Gamma=w \Gamma$ of the action of $g$ on
$G/\Gamma$. That is $w^{-1}gw\in \Gamma$. We obtain:

$$\sum_{g} \phi(g) ({\rm Trace\,} \pi)(g)
=\sum_{g}\vert\{w\in G/\Gamma; w g w^{-1}\in  \Gamma\}\vert
\phi(g).$$ Writing $g=w \gamma w^{-1}$ for some $\gamma$ in
$\Gamma$ and changing variables, we obtain:
$$<{\rm Trace\,}(Ind_{\Gamma}^G1),\phi>=\sum_{\gamma\in \Gamma/\sim}c(\gamma)<O_{\gamma},\phi>$$
where the equivalence sign denote  the action of $\Gamma$ on
$\gamma$ by conjugation, and the constant
$c(\gamma)=\frac{|G(\gamma)|}{|\Gamma(\gamma)|}$ is the quotient
of the cardinal of the stabilizers of $\gamma$ in $G$ and in
$\Gamma$.

If now we decompose the quasi-regular representation
$$Ind_\Gamma^G1=\oplus m_i \pi_i$$ in sum of irreducible representations $\pi_i$ with
multiplicities $m_i$, we obtain
$$\sum_i m_i <{\rm Trace}(\pi_i),\phi>=\sum_{\gamma\in \Gamma/\sim}
c(\gamma)<O_\gamma,\phi>.$$

This is the baby version of the trace formula.

\subsubsection{ Teenager version of the trace formula}

The preceding definitions make sense for a locally compact group
$G$ with Haar measure $dg$.

\begin{definition}\label{Tphi}
If $\phi$ is a function on $G$, and $\pi$ a representation of $G$
in a Hilbert space $V$,  the operator $\pi(\phi)$ is the operator
(provided integrals are convergent)
$$\pi(\phi)=\int_G \phi(g)\pi(g)dg.$$

\end{definition}
A representation $\pi$ in a Hilbert space $V$ may have a
distributional trace:
$$<{\rm Trace\,}\pi,\phi>={\rm Trace}_V(\pi(\phi))={\rm Trace}_V  (\int_G \phi(g)\pi(g)dg).$$

\begin{definition}

Let $O_\gamma$ be a closed orbit in $G$. Orbital integrals may
also (in good cases) be defined as distributions:
$$<O_\gamma,\phi>=\int_{O_\gamma}\phi(u)d_\gamma u$$
with respect to a invariant measure $d_\gamma u$ on $O_\gamma$.
\end{definition}

 More generally, for a packet $\{\gamma\}$ of conjugacy classes,
$<O_{\{\gamma\}},\phi>$ is a sum (eventually weighted) sum of
orbital integrals in the packet.

Thus we have two important sets of invariant distributions on the
locally compact group $G$: the characters and the orbital
integrals. One set is on the side of ``representation theory", the
other one is on the side of ``geometry". Harmonic analysis on $G$
mainly consists on understanding how to write an invariant
distribution in one set in function of the other set.

Let $\Gamma$ be a discrete subgroup with compact quotient
$G/\Gamma$, then the representation $\pi:=Ind_{\Gamma}^G1$ has a
trace and
$$<\Tr \,\pi,\phi>=\sum_{\gamma\in \Gamma/\sim} c(\gamma)<O_\gamma,\phi>$$
leading to the relation

\begin{equation}\label{traceformula}
\sum_i m_i <\Tr\, \pi_i,\phi>=\sum_{\gamma\in \Gamma/\sim}
c(\gamma)<O_\gamma,\phi>,
\end{equation}
if $\pi=\oplus m_i\pi_i$.
 The left hand side is the representation
side, while the left hand side is the geometric side.

{\bf Example: $\R/\Z$.}

{\it On the representation side, we decompose $L^2(\R/\Z)$ in a
discrete sum of representations using Fourier series, while on the
geometric side we just take the value of $\phi$ at integers. Thus
we obtain $$\sum_{n\in \Z}\int_{\R} e^{2i\pi n
x}\phi(x)dx=\sum_{n\in \Z}\phi(n).$$ This is Poisson formula.} We
have to make some reasonable assumptions on $\phi$ for the Poisson
formula to hold: for example, if $\phi$ in the Schwartz space,
this is of course true.

Arthur-Selberg trace formula is the generalization of the simple
formula (\ref{traceformula}) to cases of non compact quotients as
$SL(2,\R)/SL(2,\Z)$. It is a main tool in automorphic forms. It
allows to deduce from geometric statements the existence of some
wanted representations. Let us very roughly explain how it is used
in lifting questions.

\subsection{Transfer and the fundamental lemma}

 Now let  $G_1$ and $G_2$ be two locally compact groups.
  Assume that there is a natural map
$N:=Conj(G_1)\to Conj(G_2)$ (or more generally a map between
packets of conjugacy classes). For example, let $G_1$ be a
subgroup of $G_2$. Then we associate to a conjugacy class $O$ of
$G_2$  the packet $O\cap G_1$ of conjugacy class of $G_1$. However
in general the groups $G_1, G_2$ do not need to be related.

{\bf Example.} {\it  Let $G_1=SU(2)$ and $G_2=SL(2,\R)$. Conjugacy
classes of $G_1$ are
classified by the matrices  in $$T:=\{g(\theta)=\left(%
\begin{array}{cc}
  cos(\theta)& sin(\theta) \\
  -sin(\theta) & cos(\theta)  \\
\end{array}\right)\}.$$

We associate to the conjugacy class of $g(\theta)$ in $SU(2)$ the
conjugacy class in $SL(2,\R)$ of the same matrix $g(\theta)$ which
is also in $SL(2,\R)$ (more precisely, the packet of the two
conjugacy classes of $g(\theta)$ and $g(-\theta)$).}

 Neither of
the groups $SL(2,\R)$ and $SU(2)$ is a subgroup of the other. In
Langlands terminology, $SL(2,\R)$ is an endoscopic group of
$SU(2)$, and they share the same Cartan subgroup $T$.

Let $G_1$, $G_2$ be reductive groups over a local field with a map
$N:Conj(G_1)\mapsto Conj(G_2)$. Let $\pi$ be an irreducible
representation of $G_1$. Then the character of $\pi$  is a
distribution, which is  is given  on a dense open set (the union
of the closed conjugacy classes of maximal dimension) by
integration against a smooth function. If $\pi$ is an irreducible
representation of $G_1$, we say that a representation $\pi'$ (or a
finite sum  of irreducible representations of $G_2$) is the lift
of $\pi$ if  the character of $\pi'$ coincides with the character
of $\pi$ on matched conjugacy classes  up to a ``transfer factor".

{\bf Example}

{\it In the example of $SU(2)$, given
 the finite dimensional irreducible representation $\pi_d$ of $SU(2)$
of dimension $d+1$ with character $${\rm Trace\,}
\pi_d(g(\theta))=\sum_{|j|\leq d; j=d \,mod 2 } e^{i j\theta}$$ we
can find a sum $\pi'_d$ of two infinite dimensional irreducible
representations   of $SL(2,\R)$ such that
 $${\rm Trace\,} \pi'_d(g(\theta))=\sum_{|j|> d; j=d \,mod 2} e^{i j\theta}$$

As distributions, these characters (up to sign) coincide on the
conjugacy classes $g(\theta),  \theta \notin \{0,\pi\}$. Indeed
(for $d$ even) $$\Tr(\pi_d)(g(\theta))+\Tr(\pi'_d)(g(\theta))=
\sum_{a\in 2\Z}
e^{ia\theta}=\delta_0(\theta)+\delta_\pi(\theta)$$.}

For reductive groups $G_1$, $G_2$ over $\Q_p$, we see that there
may be  a relation between conjugacy classes if $G_1$ and $G_2$
share a common subgroup $T$. More exactly, let $T$ be an abelian
group with two homomorphisms in $G_1$, $G_2$ with images  a
maximal abelian subgroup of semi-simple elements of $G_1$
(respectively $G_2$) (Cartan subgroups of $G_1$, $G_2$). Then an
element $\gamma$ of $T$ give rise to conjugacy class
$O_\gamma^{G_1}$ and $O_\gamma^{G_2}$ in $G_1$ and in $G_2$ (or
packets). The transfer conjecture says (roughly) that if $f_1$ is
a function on $G_1$, there exists a function $f_2$ on $G_2$ such
that  the orbital integrals of $f_1$ on $O_\gamma^{G_1}$ coincide
with the orbital integral of $f_2$ on $O_\gamma^{G_2}$ (some
factors are needed) for all $\gamma\in T$. One particular case is
as follows.
 Assume
$K_1=G_1(\Z_p)$, $K_2=G_2(\Z_p)$ are maximal compact subgroups of
$G_1$, $G_2$. Then one hope that the orbital integrals of the
characteristic functions $1_{K_1}$ and $1_{K_2}$ are related on
$O_\gamma^{G_1}$ and $O_\gamma^{G_2}$. In short $1_{K_1}$ is the
transfer of $1_{K_2}$ (and Waldspurger proved that transfer of
other functions follows automatically). This is the fundamental
Lemma: (still a conjecture)

{\bf Fundamental Lemma} (Very roughly)

$$<O_\gamma^{G_1},1_{K_1}>=\mbox{(transfer constant)}*<O_\gamma^{G_2},1_{K_2}>.$$

  The precise statements have been formulated by
Langlands-Shelstad. This ``matching of orbital integrals" is one
of the tools to obtain liftings of automorphic forms. Indeed this
equality (introduced by  Labesse-Langlands) allows to compare the
contributions of spherical functions in the trace formula at
almost all places.

We give some explanations only in the case of  linear groups and
unitary groups, and we give a simplified version.

\subsubsection{Fundamental lemma and counting lattices}

I will  state the fundamental lemma for the example of the linear
group. It is equivalent to a problem of counting lattices stable
under a transformation $\gamma$. In the next subsection, I will
give
 the calculation of Labesse-Langlands in a case where it can
be done ``by hand".

Let $F$ be a local field, and $E=\oplus_{i=1}^n F e_i$ be the
standard $n$ dimensional vector space over $F$, with standard
$O_F$-lattice $L_0=\oplus_{i=1}^n O_F e_i$. Then the set
$GL(n,F)/GL(n,O_F)$ is in bijection with the set $\CL$ of lattices
$L$ over $O_F$: $L=\oplus_{i=1}^n O_F A_i$, with $A_i$ independent
vectors in $E$.  We denote  $GL(n,O_F)$ by $K$ as usual.
 We further
decompose $\CL$ in the unions of $\CL(r)$ where $r$ is an integer
(mod $n$) and

$$\CL(r)=\{L\in \CL; \,\,{\rm length\,} (L/L\cap
L_0) -{\rm length\,} (L_0/L\cap L_0)=r \,{\rm mod\,}n \}.$$

{\bf Example:\,}{\it Let us consider $E=\Q_p e_1\oplus\Q_p e_2$.
If
 $$L(x,y)=\Z_p e_1\oplus \Z_p(xe_1+y e_2),$$ then we obtain all lattices up to
homotheties ($L\mapsto p^h L$), when  $y$ varies in
$\Q_p^*/\Z_p^*$, and $x$ in $ \Q_p/\Z_p$. If $val(y)$ is even,
then $L(x,y)$ is in $\CL(0)$, while if $val(y)$ is odd, $L(x,y)$
is in $\CL(1)$. }

\bigskip

Let $\gamma\in GL(n,F)$ be an element such that its characteristic
polynomial is irreducible. Let  $\CL(\gamma)$ be the  set of
lattices stable by $\gamma$
$$\CL(\gamma):=\{L \in \CL; \gamma(L)=L\}.$$

Then this number, modulo homothetic lattices, is finite and (for
an adequate measure normalization), the orbital integral just
counts the number of elements in $\CL(\gamma)$, modulo
homotheties.
$$<O_\gamma, 1_K>={\rm cardinal\,}(\CL(\gamma)/Z),$$
where $Z$ denotes  the action of  $\Z$ by homotheties $L\mapsto
p^h L$.

{\bf Example:\,}{\it
 Let $E=Q_p e_1\oplus \Q_p e_2$ (with
$p\neq 2$). Let $\gamma\in GL(2,\Z_p)$ of the form
$$\gamma:=\left(%
\begin{array}{cc}
  a & b\delta \\
  b & a \\
\end{array}%
\right)$$ where  $\delta\in \Z_p^*$ is not a square, $a\in \Z_p^*$
and $val(b)>0$. Then the characteristic polynomial of $\gamma$ is
irreducible (here $b\neq 0$). It is easy to see that $L(x,y)$ is
stable by $\gamma$, if and only if
 \begin{equation}
 -val(b)\leq val(y)\leq val(b)
\end{equation}
\begin{equation}
  \frac{1}{2}(val(y)-val(b))\leq val(x)
\end{equation}

As $y$ varies in $\Q_p^*/\Z_p^*\sim \Z$, and $x$ in $\Q_p/\Z_p$,
we see that the set $\CL(\gamma)$ modulo homotheties  is finite.}

To formulate the fundamental lemma, we need to introduce twisted
integrals. Instead of giving the integral definition, we just give
what  it computes:

\begin{lemma}
Let $\kappa$ an integer $0\leq \kappa\leq (n-1)$. The twisted
orbital integral $<O_{\gamma,\kappa},1_K>$ is equal to
$$<O_{\gamma,\kappa},1_K>=\sum_{r=0}^{n-1}e^{2i\pi r\kappa/n}{\rm
cardinal}((\CL(\gamma)\cap \CL(r))/Z).$$

If $\kappa=0$, we obtain just the usual orbital integral
$<O_{\gamma,\kappa},1_K>$.
\end{lemma}

Let $F'$ be an extension of $F$ of degree $d$, where $n=dm$. It is
clear that an element of $GL(m,F')$ gives rise to an element of
$GL(n,F)$  (a vector space of dimension $m$ over $F'$ is of
dimension $n=md$ over $F$). Let $\xi:GL(m,F')\to GL(n,F)$ be the
corresponding homomorphism. Let $\gamma'\in GL(m,F')$ such that
$\gamma=\xi(\gamma')$ has an irreducible characteristic
polynomial. We denote by $G_1=GL(m,F')$ and by $G_2=GL(n,F)$, the
compact groups $K_1$ and $K_2$ being as usual. Then the
fundamental lemma asserts

{\bf Fundamental Lemma for linear groups}
$$<O_{\gamma'}^{G_1},1_{K_1}>=\Delta(\gamma)<O_{\gamma,m}^{G_2},1_{K_2}>.$$

Here $\gamma=\xi(\gamma')$ and is assumed to have an irreducible
characteristic polynomial. The factor $\Delta(\gamma)$ is an
explicit function of the eigenvalues of $\gamma$ and is called the
transfer factor.

This ``lemma" has been proved by Waldspurger \cite{walds}.

\subsubsection{Labesse-Langlands simplest example}

The following case is  the simplest case of the cases encountered
by Labesse-Langlands. We will see that the computation is possible
to do by hand, but relies already on wonderful  cancellations.

Let $n=2,d=2,m=1$, $G_1=GL(2,\Q_p)$ and $G_2=GL(1,F')$, where $F'$
is a quadratic extension of $\Q_p$.

Thus consider $F'=\Q_p(\sqrt{\delta})$ where $\delta\in \Z_p^*$ is
not a square. Then for $\gamma'=a+b \sqrt{\delta}$ non zero in
$F'$, the element $\xi(\gamma')\in GL(2,\Q_p)$ is equal to
$\gamma:=\left(%
\begin{array}{cc}
  a & b\delta \\
  b & a \\
\end{array}%
\right)$.

Then the fundamental lemma asserts that
$$<O_{\gamma,1}^{G_1},1_{K_1}>=(-p)^{val(b)} <1_{O_{F'}^*},\gamma'>.$$

Let us prove the fundamental lemma above: The first member is $1$
or $0$ according to the fact that $\gamma'\in O_{F'}^*$ or not. It
is easy to see that if $\gamma'$ is not in $O_{F'}^*$, there is no
invariant lattice under $\gamma$, so that the formula is true in
this case. Assume now that $\gamma$ is in $O_{F'}^*$ and that
$\gamma=a+b\sqrt{\delta}$ with $a,b\in \Z_p$, $val(b)>0$ and
$val(a)=0$. From the preceding description of all the lattices
$L(x,y)$ stable by $\gamma$, we have first to count for $val(y)$
fixed, the number $q(y)$ of $x\in \Q_p/\Z_p$ with $val(x)\geq
\frac{1}{2}(val(y)-val(b))$.

 As $val(x)$ is an integer, this is $p^{s}$ where $s$ is
 $\frac{1}{2}(val(y)-val(b))$ if $val(y)$ and $val(b)$ have same parity, or
 $\frac{1}{2}(val(y)+1-val(b))$ if $val(y)$ and $val(b)$ do not
 have the same parity.
 When we
vary $val(y)$  from $val(b)$ to $-val(b)$, this number $q(y)$
remains  the same for $2$ consecutive values of $val(y)$.

 As we
compute twisted orbitals, we take the alternate sum of the numbers
$q(y)$. It follows that all terms cancel except for the last
possible value of $val(y)=-val(b)$, and the formula follows.

\subsubsection{Unitary groups}

Assume now $F'$ is a non ramified quadratic extension of $F$, and
let $E'=\oplus_{i=1}^nF' e_i$. Let $J$ be a matrix   written as
$$J=\{\left(\begin{array}{ccccc}
  c_1& 0&0&0&0 \\
  0&c_2 & 0&0&0  \\
0&0&*&0&0\\
0&0&0&*&0\\
0&0&0&0&c_n\\
\end{array}\right)\}$$
where $c_i\in O_F^*$.
 Denote $x\mapsto
\overline{x}$ the conjugation in $F'$. Consider the unitary group
$$G=U(n,J,F)=\{M\in GL(n,F'); MJ\overline{M}^t=J\}$$
and let $K=GL(n,O_{F'})\cap G.$ The matrix $J$ determines an
hermitian form $q_J$ on $E'$, and $U(n,J,F)$ leaves this form
invariant.

 Take
$$\gamma=\{\left(\begin{array}{ccccc}
  \gamma_1& 0&0&0&0 \\
  0&\gamma_2 & 0&0&0  \\
0&0&*&0&0\\
0&0&0&*&0\\
0&0&0&0&\gamma_n\\
\end{array}\right)\}$$
where $\gamma_i\in F'$ satisfies $\gamma_i\overline{\gamma_i}=1$.

Then the orbital integral $<O_{\gamma}^G,1_K>$  is
 the number of $O_{F'}$-lattices $L\subset E'$ stable by $\gamma$
 and self-dual with respect to the quadratic form $q_J$.

Let us consider the packet $\{\gamma\}$ of conjugacy classes of
elements $\gamma'$ of $U(n,J,F)$ conjugated to $\gamma$ in
$GL(n,F')$. There is a similar twist factor $\kappa$ and  we can
consider the corresponding sum of orbital integrals is denoted by
$<O_{\gamma,\kappa}^G,\phi>$. The sum (without weights $\kappa$)
over $\{\gamma\}$ is called the stable orbital integral and
denoted by $<SO_{\gamma}^G,\phi>$.

 Let $n=n_1+n_2$, and $E'=E'_1\oplus E'_2$.
  Let $G_1=U(n,F,J)$, $K_1=G_1\cap GL(n,O_{F'})$.
Let $G_2=U(n_1,n_2,F,J)= G_1\cap (GL(n_1,F')\times GL(n_2,F'))$
and $K_2=K_1\cap G_2$.

 Ngo Bao Chau and Laumon proved  the fundamental lemma for stable
 orbital integrals: that is
 the stable orbital integral $<SO_{\gamma}^{G_2},1_{K_2}>$
 is compared to the
similar twisted packet $<O_{\gamma,\kappa}^{G_2},1_{K_2}>.$ Most
remarkably, the proof uses equivariant cohomology of Hitchin's
moduli space.

The proof of the fundamental lemma should lead to advances on
$L$-functions attached to Shimura varieties (higher cases of the
modular curves), etc...

\section{Geometric Langlands correspondence}

A geometric analogue of the Langlands program has been formulated
by Drinfeld. Here a finite extension $E$ of $F_q(t)$ is  replaced
by the field of meromorphic functions $\CM(C)$ over a complex
compact curve $C$.

Let us first explain the analogy of the geometric program with the
number field program.

Let $P_1(C)=\{[z_1,z_2]\}/\C^*$. The ring of rational functions on
$P_1(\C)$ is $\C(T)$ with $T=z_1/z_2$. A finite extension  of
$\C(T)$ is the same thing as a covering of $P_1(\C)$. Namely if
$m:C\to P_1(\C)$ is a finite cover, we embed $\C(T)$ into
$\CM(C)$, the function field over $C$, by the map $m^*$. Thus, the
analogue of a number field will be the field $\CM(C)$ of functions
on a complex curve $C$. A place of $\CM(C)$ is simply a point $a$
of $C$, with corresponding local field $F_a$ Laurent series in
$(z-a)$ (Here $z$ is a local coordinate around $a$). We see that
this local field is very analogous to $\Q_p$ which are ``Laurent
series" in $p$.

{\bf Example} {\it If $p_1,p_2,p_3$ are distinct primes, then the
analogue of $F:=\Q(\sqrt{p_1p_2p_3})$ will be the field  of
functions on the elliptic curve
$C:=\{y^2=(z-a_1)(z-a_2)(z-a_3)\}$. Remark that the cover
$(y,z)\mapsto z$ has two points, except at $a_1,a_2,a_3$. The only
points where this cover is ramified  are the points $a_1,a_2,a_3$.
Similarly primes $p$ different from $p_1,p_2,p_3$ (and $2$) are
unramified in $\Q(\sqrt{p_1p_2p_3})$.}

\bigskip

So let $C$ be a complex curve and $\CM(C)$ its function field. We
will consider only the ``unramified" case:  the fundamental group
$G:=\pi(C)$ of the curve $C$ is a quotient of the absolute Galois
group  and we will consider only representations of this quotient
group.
 This corresponds  to the Galois group of
unramified covers of $C$. Indeed let $Y$ be an unramified cover of
$C$, then $\CM(Y)$ is an extension of $\CM(C)$. A deck
transformation $Y\to Y$ induces a transformation of $\CM(Y)$ which
is the identity on $\CM(C)$, thus the Galois group
$Galois(\CM(Y)|\CM(C))$ is the group of deck transformations. Now,
as $Y$ is unramified,  an element of $g$ induces a deck
transformation on $Y$. Thus, we have a morphism from $G$ to Galois
groups $G(\CM(Y)/\CM(C))$  for the unramified cover $Y\to C$. In
the following, only representations of $G=\pi(C)$ will be
considered, and this corresponds (in the dictionary below) to
automorphic representations everywhere non ramified (thus there
exists a canonical vector in this representation).

Here is the dictionary of analogies for this situation (the
unramified case).

\begin{itemize}

\item  {\bf 1}

\framebox{(quotient of the) Galois group of a number field $F$}
\\$\mapsto$\\ \framebox{Fundamental group of $C$}
\\

\item {\bf 2}

\framebox{Representation of the Galois group} \\$\mapsto$\\
\framebox{Local system on $C$}
\\
\item  {\bf 3}

\framebox{The double coset defined by $g=(g_v)_{v\in\CV}$ in
$GL(n,O_F)\backslash GL(n,\A_F)/GL(n,F)$}
\\$\mapsto$\\ \framebox{A vector bundle of rank $n$ on
$C$.}\\

\item  {\bf 4}

\framebox{$GL(n,O_F)\backslash GL(n,\A_F)/GL(n,F)$}$\\\mapsto$\\
 \framebox{the moduli space $Bun_n$ of rank $n$ vector bundles over $C$}.\\

\item  {\bf 5}

\framebox{A (nice) function $f$ on $GL(n,O_F)\backslash
GL(n,\A_F)/GL(n,F)$}$\\\mapsto$\\
\framebox{A perverse sheaf $\CF$ on $Bun_n$.}\\

\item{\bf 6}

\framebox{An eigenfunction of the Hecke operators}\\$\mapsto$\\
\framebox{A Hecke eigensheaf $\CF$ on $Bun_n$.}\\
\end{itemize}

Let me explain what I understood of this dictionary.

\begin{itemize}

\item

Point {\bf 2}.
 A
representation of the fundamental  group $\sigma:G\mapsto
GL(n,\C)$ leads to a flat vector bundle over $C$. This is called a
local system $E$ over $C$. Thus we know very well the geometric
analogue of representations of the Galois group. These are local
systems $E$ on $C$.\\

\item

Point {\bf 3} and {\bf 4}.
 At a place $a$ of $C$, we consider a coordinate $z$,  and we
identify $GL(n, \CM(C)_a)$ with  $GL(n,\C((z)))$. Let $V\to C$  be
an holomorphic vector bundle of rank $n$ and choose
$S_1,S_2,\ldots,S_n$ meromorphic sections of $V$ generically
independent. At each point $a$ in $C$, we take a local
trivialisation of $V_a=\C^n$ via holomorphic sections
$s_1,s_2,\ldots,s_n$ (defined  at $a$). We obtain an element
$g_a(z)\in GL(n, \CM(C)_a)$ by writing $s_i(z)=\sum_{j}
g_a^{i,j}(z)S_j(z)$. Thus an holomorphic vector bundle $V$ of rank
$n$  gives, via its transition functions, an element
$g=(g_a(z))_{a\in C}$ with $g_a(z)\in GL(n,\C((z)))$. Thus the
analogue of the space $GL(n,O_F)\backslash GL(n,\A_F)/GL(n,F)$ is
the space $Bun_n$ of equivalence classes of holomorphic vector
bundles over $C$ of rank $n$. The space $Bun_n$ is a ``stack", but
I will employ anyway terms of ``varieties". Anyway, all this is
very vague for me. The space $Bun_n$ is not connected. It is the
union of connected components
$Bun_n^d$ of vector bundles of degree $d$.\\

\item Points {\bf 5} and {\bf 6}.

Now, it is more difficult (for me) to understand what is the
analogue of the ``automorphic side", that is what is the analogue
of an automorphic  representation.  Recall that in the case where
$\pi$ was an irreducible automorphic representation, we could
singled out in $\pi$ a particular function $f=\otimes_v f_v$ in
$\pi\subset L^2(GL(n,A_F)/GL(n,F))$, at least at all unramified
places, by saying that $f_p$ was the fixed vector under $K_p$
(normalized to be $1$ at $1$). This vector $f_p$ was automatically
an eigenvector for the Hecke operators.

Now, in the geometric context, the analogue of a function on a
space $X$ is a perverse sheaf on $X$. Indeed consider the case
where $X$ is finite. A sheaf $\CF$ of $X$ is just a collection of
vector space $\CF_x$. We consider sheaves up to isomorphisms. Thus
$\CF$ is completely determined by the function $x\mapsto
\dim(\CF_x)$.
 Remark that this
operation commutes with the $6$ operations, in the Grothendieck
dictionary: If we consider the tensor product $\CF\otimes \CG$ of
``sheaves", then the corresponding function  is $f(x)g(x)$.  If
$X\subset Y$, we extend $\CF$ on $Y$ by $\CF_y=0$ if $y\notin X$,
this corresponds to extending $f$ by $0$. Similarly for
pushforward $X\mapsto Y$, the function $\pi_*(f)(y)=\sum_{x;
\pi(x)=y} f(x)$ commutes with pushforward of sheaves
$(\pi_*\CF)_y=\oplus_{x; \pi(x)=y} \CF_x$. When $X$ is a variety
defined over a finite field $F_q$, and $\CF$ a sheaf on $X$, and
$x$ a point of $X(F_q)$, the Frobenius element $Fr_q$ acts on the
cohomology groups $H_i(x,\CF)$, and we obtain a function on
$X(F_q)$, by taking the trace of the action of the Frobenius.

If $X$ is a complex variety, an irreducible perverse sheaf is
completely determined by a locally closed complex subvariety $Y$
of $X$, and a local system on $Y$, and operations on perverse
sheaves of extensions, pushforward, etc... are well defined.

Finally,  we  discuss the notion of Hecke eigensheaves.

For automorphic forms, at each prime $p$, consider the subspace
$D^i_p$ of $GL(n,\Q_p)/GL(n,\Z_p)\times GL(n,\Q_p)/GL(n,\Z_p)$
consisting of elements $\{g_1,g_2\}$ with $g_1=h_p^ig_2$ where

$$h_p^i=\left(%
\begin{array}{ccccc}
  p & 0 & 0 & 0 & 0 \\
  0 & p & 0 & 0 & 0 \\
  0 & 0 & p & 0 & 0 \\
  0 & 0 & 0 & * & 0 \\
  0 & 0 & 0 & 0 & 1 \\
\end{array}%
\right).$$

Now at each place $a$ of $C$, we consider the following subset
$$Hecke^i_a\subset Bun_n\times Bun_n$$ of vector bundles $(V,V')$
with a map $V\to V'$ isomorphism, except at the point $a$ of $C$,
where it is locally given the matrix (analogue of the Hecke
matrix)

$$\left(\begin{array}{ccccc}
  (z-a) & 0 & 0 & 0 & 0 \\
  0 & (z-a) & 0 & 0 & 0 \\
  0 & 0 & (z-a) & 0 & 0 \\
  0 & 0 & 0 & * & 0 \\
  0 & 0 & 0 & 0 & 1 \\
\end{array}%
\right).$$

We now vary  the base point $a$, and obtain a subvariety

$$H^i\subset C\times Bun_n\times Bun_n$$

Let us denote by $p_1(a,V,V')=V'$ and $p_2(a,V,V')=(a,V)$. Thus
$H^i$ operates on perverse sheafs on $Bun_n$ by

\begin{definition}
$$H^i(\CF)=(p_2)_*p_1^*(\CF).$$
\end{definition}

Thus from a perverse sheaf on $Bun_n$, we obtain a perverse sheaf
on $C\times Bun_n$.
\\
\end{itemize}

We are now ready to state what should be the analogue of the
Langlands correspondence.

Let $E$ be a local system on $C$ (a representation of the Galois
group) of rank $n$.
 The classical  Hecke operators were operators on functions. The equation
 $H_pf=a_pf$ saying that $f$ is an eigenfunction of $H_p$
 is translated by the equation on sheaves:
 $H^1\CF=E\otimes \CF$.

\begin{theorem}
For each irreducible local system $E$ on $C$ of rank $n$, there
exists an irreducible (on each component $Bun_n^d$) perverse sheaf
$Aut_E$ (the automorphic sheaf attached to $E$) on $Bun_n$ which
is a Hecke eigensheaf:
$$H^i(Aut_E)=\wedge^i(E)\otimes Aut_E.$$
\end{theorem}

When $n=1$, the theory is very simple and due to Rosenlicht, Lang
and Serre. Indeed  $Bun_n$, for $n=1$,  is just the Picard group
$Pic(C)$ of line bundles on $C$. One of the connected component of
$Pic(C)$ is the Jacobian variety $Jac(C)$ of $C$. If $g$ is the
genus of $C$, then the  fundamental group of the Jacobian of $C$
(a complex torus of dimension $g$) is isomorphic to
$H_1(C,\Z)=\pi(C)/{\rm commutators}$. Thus there is an equivalence
between local systems of rank $1$ on $C$ or on $Jac(C)$. Thus  a
local system on $E$ gives rise to a sheaf $Aut_E$ on $Pic_d(C)$.
It is possible here  to describe how $Aut_E$ looks like. Consider
the restriction of $Aut_E$ on $Pic_d(X)$ the space of line bundles
of degree $d$. If $\CL$ is a line bundle on $C$ of degree $d$, an
section $s\in H^0(C,\CL)$ gives us  $d$ points on $C$ where $s$
vanishes. Thus we obtain a map $p$ from the fiber bundle over
$Pic(C)$, with fiber at $\CL$ the projective space $[H^0(C,\CL)]$
to  $(C\times C\times \cdots \times C)/Permutations$. The sheaf
$Aut_E$ is the unique local system on $Pic(C)$ such that its pull
back is $\otimes^d E$ on $C^d/permutations$. Its existence is
deduced from the fact that fibers (projective spaces) are simply
connected.

The construction of the automorphic  bundle $Aut_E$  for a local
system $E$ of rank $n$ requires much more subtle arguments. It
takes its grounds in works of Laumon, Frenkel-Gaitsgory-Vilonen
with the final step established by Gaitsgory.

\end{document}